\documentclass[a4paper,11pt]{amsart} 
\usepackage{fourier} 
\usepackage{amssymb}
\usepackage{amsmath}
\usepackage{amsthm}
\usepackage{hyperref}
\usepackage{fullpage}
\usepackage{xcolor} 
\usepackage{listings}
\usepackage{biblatex}
\renewbibmacro{in:}{} 

\usepackage{tikz}
\usetikzlibrary{trees,arrows,automata,positioning}
\tikzstyle{every picture}=[scale=.5,inner sep=10]

\usepackage{booktabs}

\theoremstyle{plain}
\newtheorem{thm}{Theorem}[section]
\newtheorem{prop}[thm]{Proposition}

\theoremstyle{definition}

\newtheorem*{notn*}{Notation}

\newtheorem{remark}[thm]{Remark}

\newtheorem{example}[thm]{Example}

\def\M{\mathcal{M}}
\def\Mbar{\overline{\mathcal{M}}}
\def\L{\mathbb{L}}

\title{
Deligne's weight spectral sequence and tautological cohomology of the moduli space of curves 
}
\author{Jonas Bergstr\"om and Thomas Wennink}
\address{Matematiska institutionen \\ Stockholms Universitet \\ 106 91 Stockholm \\ Sweden}
\email{jonasb@math.su.se, tomwennink42@gmail.com}


\begin{document}

\begin{abstract} 
We have written a computer program that implements Deligne's pullback and pushforward weight spectral sequences to compute the weight graded pieces of the rational cohomology of moduli spaces of pointed smooth curves (as well as curves of compact type and curves with rational tails) in cases where the cohomology groups appearing in the boundary stratification of the Deligne-Mumford compactification are generated by tautological classes (and when Pixton's relations are all relations). 
The weight graded pieces are computed together with the induced action of the symmetric group permuting the points on the curves. 
Using the computer program we have determined this information in the case of genus five as well as in the case of genus three with three marked points. 
\end{abstract}

\maketitle

\section{Introduction} \label{sec:introduction}
For any $g,n \geq 0$, with $2g-2+n > 0$, let $\M_{g,n}$ and $\Mbar_{g,n}$ denote the moduli spaces of smooth, respectively stable, $n$-pointed curves of genus $g$. In this article we implement Deligne's pullback and pushforward weight spectral sequences to compute weight graded pieces of the compactly supported cohomology of $\M_{g,n}$ in cases where the $E_1$-page, which is built out of pieces of cohomology groups of $\Mbar_{g,n}$ and of strata in the boundary $\Mbar_{g,n} \setminus \M_{g,n}$, consists only of tautological cohomology. Using this implementation we proved Theorems~\ref{thm:main1} and \ref{thm:main2} below. 
In the whole of the article, cohomology will mean either $\ell$-adic \'etale cohomology or of rational Hodge structures. 

The symmetric group $\mathbb S_n$ acts on $\M_{g,n}$ and $\Mbar_{g,n}$, and it induces an action on the cohomology groups of these spaces. 
This action also respects the weight filtration $W_*$ (induced by Deligne's weight spectral sequences).  
To each partition $\lambda$ of $n$ we associate, in the usual way, an irreducible representation $s_{\lambda}$ of $\mathbb S_n$. 
What we aim to compute is then, for each $q,r$, the following decomposition of the weight graded pieces, 
\[
\mathrm{gr}_q H^r(\M_{g,n})=W_q H^r(\M_{g,n})/W_{q-1}H^r(\M_{g,n})=\bigoplus_{\lambda \vdash n} \left(U^{q,r,\lambda}_{g,n} \boxtimes s_{\lambda} \right), 
\]
for some (pure) rational Hodge structures (respectively $\ell$-adic Galois representations) $U^{q,r,\lambda}_{g,n}$. Using Poincar\'e duality we then get that
\[\mathrm{gr}_q H_c^r(\M_{g,n}) = (\mathrm{gr}_{2(3g-3+n)-q} H^{2(3g-3+n)-r}(\M_{g,n}))^{\vee} \, \L^{3g-3+n}.
\]

In general, we know that for $g\geq 1$, $H_c^k(\M_{g,n})$ vanishes for $k<2g$ when $n=0,1$, for $k<2g+1$ when $n=2$ and for $k<2g-2+n$ when $n \geq 3$, see \cite{harer}, \cite{CFP}, \cite{MSS} and \cite[Theorem 4.1]{Wong}. The results of \cite{CGP} and \cite{paynewillwacherweight2} have made the weight $0$ and $2$ part of the compactly supported cohomology of $\M_{g,n}$ amenable to computation.
Furthermore, the weight graded pieces for odd weights up to $9$ all vanish, see \cite{arbarellocornalba} and \cite[Theorem 1.1]{BFP}. 
Finally, the isomorphism in \cite[Proposition 1.7]{paynewillwacherweight11} gives a way to compute the weight~$11$ contribution to any compactly supported cohomology group of $\M_{g,n}$.

The weight graded pieces (or indeed the Betti numbers) of the cohomology of $\M_{g,n}$ have only been fully determined for a few $g,n$, when $g>0$. 
For results in $g=2$ and $n \leq 5$ see \cite{tommasi2n}, in $g=3$ and $n\leq 1$ see \cite{looijenga}, in $g=3$ and $n=2$ see \cite[Theorem 1.1]{tommasi32}, in $g=4$ and $n=0$ see \cite[Theorem 1.4]{tommasi4} and finally for $g=4$ and $n=1$ see \cite[Theorem 1.2]{zhengwong}. As far as we know, there is for $g=1$ a lack of published results. Combining the methods in \cite{peterseng1,totaro,gorinov}, Louis Hainaut has computed all results for $g=1$ and $n\leq 6$ (personal communication). All these results (except for the case of $g=2$ and $n=5$), have been verified using our implementation and are listed in Section~\ref{sec:hodgeserre}.

By \cite[Theorem~1.5]{canninglarsonpaynewillwacherPol}, if $3g+2n<25$ (a range that covers all results mentioned above as well as in the following two theorems), the cohomology will be of Tate type. This means that the odd weight graded pieces vanish and that a $2q$-weight graded piece of a cohomology group will be a direct sum of $q$th tensor products of the Tate class $\L=H^2(\mathbb P^1)$, denoted~$\L^q$. The associated pure rational Hodge structure of $\L$ has Hodge weight $(1,1)$ and the associated $\ell$-adic Galois representation is the cyclotomic character. 

In the following two theorems the weight filtration of each cohomology group is concentrated in a single weight except for $H^9(\M_{3,3})$.

\begin{thm} \label{thm:main1} 
The cohomology of $\M_{5}$ is of Tate type and equals,
\begin{multline*}
H^0(\M_{5})=1, \quad 
H^2(\M_{5})=\L,  \quad 
H^4(\M_{5})=\L^2,  \quad 
H^5(\M_{5})=\L^3,\\
H^6(\M_{5})=\L^3, \quad 
H^7(\M_{5})=\L^4,  \quad 
H^{14}(\M_{5})=\L^{12}.    
\end{multline*}
\end{thm}
\begin{remark} The rational Chow ring of $\M_{5}$ equals $\mathbb Q[\kappa_1]/\kappa_1^4$, see \cite{izadi}. The cycle class map from the rational Chow ring into rational cohomology is in this case injective, giving the contributions $1$, $\L$, $\L^2$ and $\L^3$ in degrees $0$, $2$, $4$ and $6$ in Theorem~\ref{thm:main1}. 

The moduli space $\M_{5}$ can be stratified into the hyperelliptic locus $\mathcal H_{5}$, the trigonal locus $\mathcal T_5$ and an affine open locus of general curves $\mathcal O_5$. The rational cohomology of  $\mathcal H_{5}$ equals that of a point. In \cite{zhengtrigonal5}, the rational cohomology of $\mathcal T_5$ was determined. Orsola Tommasi and Angelina Zheng then conjectured (building upon earlier computations of Orsola Tommasi using the Gorinov-Vassiliev method applied to $\mathcal O_5$), that the rational cohomology of $\mathcal O_5$ would consist of a class $1$ in degree $0$ and $\L^3$ in degree $5$. Another way to prove Theorem~\ref{thm:main1} would be to prove this conjecture about the rational cohomology of $\mathcal O_5$. 
\end{remark}

\begin{thm} \label{thm:main2} 
The cohomology of $\M_{3,3}$ is of Tate type and equals,
\begin{multline*}
H^0(\M_{3,3})=s_3, \quad 
H^2(\M_{3,3})=(2s_3+s_{2,1})\L, \quad 
H^4(\M_{3,3})=(s_3+s_{2,1})\L^2, \quad
H^5(\M_{3,3})=2(s_3+s_{2,1})\L^3,  \\
H^6(\M_{3,3})=s_3\L^6, \quad 
H^8(\M_{3,3})=2(s_3+s_{2,1})\L^7, \quad 
\mathrm{gr}_{16} H^9(\M_{3,3})=2s_3\L^8, \quad
\mathrm{gr}_{18} H^9(\M_{3,3})=s_3 \L^9. 
\end{multline*}
\end{thm}

\begin{remark}
The $\mathbb{S}_n$-equivariant Euler characteristic of $\M_{5}$ and $\M_{3,3}$ in the Grothendieck group of mixed rational Hodge structures (or $\ell$-adic Galois representations) had already been computed and will appear in forthcoming work by the first author together with Samir Canning, Dan Petersen and Johannes Schmitt. 
\end{remark}

Our program builds upon the SageMath module {\tt admcycles}, see \cite{admcycles}, and improves parts of it. Our code is available at \cite{DeligneWeightSpectral} and we plan to integrate it into the {\tt admcycles} package.

\begin{remark}
Using {\tt admcycles} (and checking that Pixton's relations are all relations, see further in sections~\ref{sec:taut} and \ref{sec:implementation}), one can show that in degrees $0$, $2$, and $4$, the cohomology of $\M_{3,3}$ is tautological.  
\end{remark}

Whenever we compute the weight graded pieces of the cohomology of $\M_{g,n}$,
we also compute the weight graded pieces of the cohomology of $\M^\mathrm{ct}_{g,n}$ and $\M^\mathrm{rt}_{g,n}$, the moduli space of curves of compact type, respectively with rational tails.
This is included in our list of results at \cite{DeligneWeightSpectral}.
The idea to apply Deligne's weight spectral sequence to the moduli space of curves of compact type comes from forthcoming work by Amy Bradford.

\subsection*{Acknowledgements}
We are most grateful to Amy Bradford, Samir Canning, Hannah Larson, Sam Payne, Dan Petersen, Johannes Schmitt, Orsola Tommasi, Thomas Willwacher and Angelina Zheng for valuable discussions and comments, and to Samir Canning for suggesting to implement Deligne's weight spectral sequences. We thank the referees for their careful reading and helpful comments.  
During this work, the second author was supported by a Sverker Lerheden fellowship for postdoctoral studies. 

\section{Deligne's pushforward weight spectral sequence} 
\label{sec:deligne}
We will follow \cite[Section 7]{arbarello2009divisors} and \cite[Section 2.3]{paynewillwacherweight2}, but write the differential in a slightly more explicit way which better corresponds to our concrete computations. 

To an $n$-pointed genus $g$ stable graph $\Gamma$, with $V_\Gamma$ and $E_\Gamma$ denoting its set of vertices respectively edges, we associate a product of moduli spaces and a gluing map
\[\Mbar_\Gamma=\prod_{v\in V_\Gamma}\Mbar_{g(v),n(v)}
,\qquad
\xi_\Gamma:\Mbar_\Gamma\rightarrow\Mbar_{g,n}.\]
For any stable graph $\Gamma$, we have an action of $\mathrm{Aut}(\Gamma)$ on $\Mbar_\Gamma$, inducing an action on $H^k(\Mbar_\Gamma)$ for each~$k$. 

Deligne's pushforward weight spectral sequence \cite[Th\'eor\`eme (3.2.5)]{deligne} gives, in the case of $\M_{g,n} \subset \Mbar_{g,n}$, an $E_1$-page of the form, 
\[E_1^{-p,q} =\bigoplus_{|E_\Gamma|=p}(H^{q-2p}(\Mbar_\Gamma) \otimes \mathrm{det}(E_\Gamma))^{\mathrm{Aut}(\Gamma)} \otimes \L^p,\]
where $\mathrm{Aut}(\Gamma)$ acts on $E_\Gamma$ by permutations.
This spectral sequence degenerates at the second page and abuts to the weight $q$ graded piece, 
\[E_2^{-p,q} =\mathrm{gr}_q H^{q-p}(\mathcal{M}_{g,n}).
\]

For every isomorphism class of stable graphs we choose a representative $\Gamma$ and we also fix an ordering of its edges $\alpha_\Gamma:E_\Gamma\rightarrow \{1,\ldots,|E_\Gamma|\}$.
For the stable graph $\Gamma$ we let $\Gamma_i$ denote the graph obtained by contracting its $i$-th edge. 
It inherits an ordering $\beta_{\Gamma_i}:E_{\Gamma_i}\rightarrow\{1,\ldots,i-1,i+1,\ldots,|E_\Gamma|\}$.
Let $\Gamma'_i$ be the chosen representative of the isomorphism class of $\Gamma_i$. Fix an isomorphism $\phi:\Gamma_i\rightarrow\Gamma'_i$.

For each of the edge contractions we have a gluing map
\[\Mbar_\Gamma\xrightarrow{\xi_{\Gamma,\Gamma_i}}\Mbar_{\Gamma_i},\]
and, for each $k$, we get the differential restricted to $\Gamma$, 
\begin{align}
    \nonumber
	d_1:\bigl(H^r(\Mbar_\Gamma) \otimes \mathrm{det}(E_\Gamma)\bigr)^{\rm{Aut}(\Gamma)}&\rightarrow 
    \bigoplus_{|E_{\Gamma'}|=|E_\Gamma|-1} \bigl(H^{r+2}(\Mbar_{\Gamma'}) \otimes \mathrm{det}(E_{\Gamma'})\bigr)^{\rm{Aut}(\Gamma')}\\
    \label{dwss_diff}
    x & \mapsto \sum_{i=1}^{|E_\Gamma|}\frac{1}{|\mathrm{Aut}(\Gamma'_i)|}\sum_{\psi\in\mathrm{Aut}(\Gamma'_i)}\textrm{sgn}(\sigma_{i,\psi})\cdot(\psi\circ\phi\circ\xi_{\Gamma,\Gamma_i})_*(x),
\end{align}
where $\sigma_{i,\psi}:=\alpha_{\Gamma'_i}\circ \psi_E\circ\phi_E\circ\beta^{-1}_{\Gamma_i}$ 
is the permutation of $\{1,\ldots,|E_\Gamma|\}$ that sends $i$ to $|E_\Gamma|$ and that acts by $\alpha_{\Gamma'_i}\circ \psi_E\circ\phi_E\circ\beta^{-1}_{\Gamma_i}$ elsewhere. Here $\psi_E,\phi_E$ denote the induced maps on edges. The contribution of $\textrm{sgn}(\sigma_{i,\psi})$ comes from the conversion of $\mathrm{det}(E_\Gamma)$ coefficients into $\mathrm{det}(E_{\Gamma'})$ coefficients. 

\begin{example}
	Consider the following stable graphs with ordered edges
	\[
		\Gamma =
	\begin{tikzpicture}[baseline={([yshift=-.8ex]current bounding box.center)},el/.style = {inner sep=2pt, align=left, sloped},every child node/.style={inner sep=1,font=\tiny}]
      \tikzstyle{level 1}=[counterclockwise from=0,level distance=9mm,sibling angle=120]
			\node (A0) [draw,circle,inner sep=1] at (0:-1.5) {$\scriptstyle{1}$};
      \tikzstyle{level 1}=[counterclockwise from=120,level distance=9mm,sibling angle=120]
      \node (A1) [draw,circle,inner sep=1] at (0:0) {$\scriptstyle{2}$};
      \tikzstyle{level 1}=[counterclockwise from=-60,level distance=9mm,sibling angle=120]
      \node (A2) [draw,circle,inner sep=1] at (0:1.5) {$\scriptstyle{0}$} child {node {1}} child {node {2}};

			\path (A0) edge [bend left=0.000000] node[el,above,font=\tiny] {$e_1$}(A1);
			\path (A1) edge [bend left=0.000000] node[el,above,font=\tiny] {$e_2$} (A2);
	\end{tikzpicture}
,\quad
\Gamma'_1 =
	\begin{tikzpicture}[baseline={([yshift=-.8ex]current bounding box.center)},el/.style = {inner sep=2pt, align=left, sloped},every child node/.style={inner sep=1,font=\tiny}]
      \tikzstyle{level 1}=[counterclockwise from=120,level distance=9mm,sibling angle=120]
      \node (A1) [draw,circle,inner sep=1] at (0:-.75) {$\scriptstyle{3}$};
      \tikzstyle{level 1}=[counterclockwise from=-60,level distance=9mm,sibling angle=120]
      \node (A2) [draw,circle,inner sep=1] at (0:.75) {$\scriptstyle{0}$} child {node {1}} child {node {2}};

			\path (A1) edge [bend left=0.000000] node[el,above,font=\tiny] {$e_1$} (A2);
	\end{tikzpicture}
	, \quad
		\Gamma'_2 =
	\begin{tikzpicture}[baseline={([yshift=-.8ex]current bounding box.center)},el/.style = {inner sep=2pt, align=left, sloped},every child node/.style={inner sep=1,font=\tiny}]
      \tikzstyle{level 1}=[counterclockwise from=0,level distance=9mm,sibling angle=120]
			\node (A0) [draw,circle,inner sep=1] at (0:-.75) {$\scriptstyle{1}$};
      \tikzstyle{level 1}=[counterclockwise from=-60,level distance=9mm,sibling angle=120]
      \node (A1) [draw,circle,inner sep=1] at (0:.75) {$\scriptstyle{2}$} child {node {1}} child {node {2}};

			\path (A0) edge [bend left=0.000000] node[el,above,font=\tiny] {$e_1$}(A1);
	\end{tikzpicture}
.\]

	On $\Gamma$, the differential map can be restricted to a map
	\[
    \bigl(H^*(\Mbar_\Gamma) \otimes \mathrm{det}(E_\Gamma)\bigr)^{\rm{Aut}(\Gamma)}
    \rightarrow
    \bigl(H^*(\Mbar_{\Gamma'_1}) \otimes \mathrm{det}(E(\Gamma'_1))\bigr)^{\rm{Aut}(\Gamma'_1)}
    \bigoplus
    \bigl(H^*(\Mbar_{\Gamma'_2}) \otimes \mathrm{det}(E(\Gamma'_2))\bigr)^{\rm{Aut}(\Gamma'_2)} 
    \]
    such that
    \[
	1\otimes 1\otimes 1 \mapsto
	\biggl( \bigl[
	\begin{tikzpicture}[baseline={([yshift=-.8ex]current bounding box.center)},el/.style = {inner sep=2pt, align=left, sloped},every child node/.style={inner sep=1,font=\tiny}]
      \tikzstyle{level 1}=[counterclockwise from=0,level distance=9mm,sibling angle=120]
			\node (A0) [draw,circle,inner sep=1] at (0:-.75) {$\scriptstyle{1}$};
      \tikzstyle{level 1}=[counterclockwise from=0,level distance=9mm,sibling angle=120]
      \node (A1) [draw,circle,inner sep=1] at (0:.75) {$\scriptstyle{2}$} child {node {1}};
			\path (A0) edge [bend left=0.000000] (A1);
	\end{tikzpicture}
	\bigr]
	\otimes 1\biggr)
	+(-1) \cdot  
	\biggl( 1 \otimes
	\bigl[
	\begin{tikzpicture}[baseline={([yshift=-.8ex]current bounding box.center)},el/.style = {inner sep=2pt, align=left, sloped},every child node/.style={inner sep=1,font=\tiny}]
      \tikzstyle{level 1}=[counterclockwise from=180,level distance=9mm,sibling angle=120]
      \node (A1) [draw,circle,inner sep=1] at (0:-.75) {$\scriptstyle{2}$} child {node {3}}; child {node {}};
      \tikzstyle{level 1}=[counterclockwise from=-60,level distance=9mm,sibling angle=120]
      \node (A2) [draw,circle,inner sep=1] at (0:.75) {$\scriptstyle{0}$} child {node {1}} child {node {2}};
			\path (A1) edge [bend left=0.000000] (A2);
	\end{tikzpicture}
	\bigr] \biggr),
    \]
	where $1$ represents the fundamental class.
\end{example}

\subsection{Compact type and rational tails} \label{sec:ct-rt}
A stable curve is of compact type if its dual graph has no cycles.
We say it is a curve with rational tails if it has one irreducible component of full genus, because then all other components are chains (tails) of rational curves.

We consider Deligne's pushforward weight spectral sequence in the cases
$\M^\mathrm{ct}_{g,n}\subset \Mbar_{g,n}$
and
$\M^\mathrm{rt}_{g,n} \subset \Mbar_{g,n}$.
The first page of the spectral sequence is similar to the one for 
$\M_{g,n} \subset \Mbar_{g,n}$.
The difference is that in the case of curves of compact type we sum over stable graphs without bridges and 
in the case of curves with rational tails we sum over stable graphs without any rational tails.

This can be proven by applying the description of how the boundary strata of the moduli space of stable curves meet, but restricted to the subset of boundary strata of the moduli space of curves of compact type (respectively rational tails), see for instance \cite[Chapter XII, \S 10]{arbarellocornalbagriffiths}. In the case of compact type, more details will be found in the forthcoming thesis of Amy Bradford.

\section{Deligne's pullback weight spectral sequence} \label{sec:gk}
When doing concrete computations, it is helpful to use both Deligne's pushforward weight spectral sequence and its Poincaré dual version, which has an $E_1$-page of the form, 
\[\tilde E_1^{p,q}=\bigoplus_{|E_\Gamma|=p} \bigl(H^{q}(\Mbar_\Gamma) \otimes \mathrm{det}(E_\Gamma)\bigr)^{\mathrm{Aut}(\Gamma)},\]
 see \cite[Example 3.5]{petersen}. This spectral sequence also degenerates at the second page and abuts to the weight $q$ graded piece, 
\[\tilde E_2^{p,q} =\mathrm{gr}_q H_c^{p+q}(\mathcal{M}_{g,n}).
\]
Fixing $q$, the complex $\tilde E_1^{p,q}$ can also be identified with the Getzler-Kapranov complex $\mathrm{GK}^q_{g,n}$. For more details see \cite[Section 2.3]{paynewillwacherweight2} and \cite{GK}. 
Similarly to Section~\ref{sec:deligne}, the differential 
\[ d_1: \tilde E_1^{p,q} \rightarrow \tilde E_1^{p+1,q} \]
restricted to a graph $\Gamma$ is given by 
\begin{equation}
\label{GK_delta}
x \mapsto \sum_{v\in V_\Gamma}\sum_{\hat\Gamma}\frac{1}{|\mathrm{Aut}(\hat\Gamma')|}\sum_{\psi\in\mathrm{Aut}(\hat\Gamma')}\textrm{sgn}(\sigma_{\hat\Gamma,\psi})\cdot(\xi_{\hat\Gamma,\Gamma}\circ\phi\circ \psi)^*(x).
\end{equation}
Here the second sum is over all stable graphs $\hat\Gamma$ obtained by either splitting the vertex $v$ or adding a loop to the vertex $v$ and lowering its genus by one.
This operation is dual to contracting an edge in the sense that by contracting this new edge of $\hat\Gamma$, we obtain $\Gamma$ again.
From the ordering of $\Gamma$ we inherit the ordering $\gamma_{\hat\Gamma}:E_{\hat\Gamma}\rightarrow\{1,\ldots,|E_\Gamma|+1\}$ obtained by sending the new edge to $|E_\Gamma|+1$.
Let $\hat\Gamma'$ be the chosen representative of the isomorphism class of $\hat\Gamma$ and fix an isomorphism $\phi:\hat\Gamma\rightarrow\hat\Gamma'$. 
Now $\sigma_{\hat\Gamma,\psi}:=\alpha_{\hat\Gamma'}\circ \psi_E\circ\phi_E\circ\gamma^{-1}_{\hat\Gamma}$.

\begin{example}
	Let	$\Gamma =
	\begin{tikzpicture}[baseline={([yshift=-.8ex]current bounding box.center)},el/.style = {inner sep=2pt, align=left, sloped},every child node/.style={inner sep=1,font=\tiny}]
      \tikzstyle{level 1}=[counterclockwise from=0,level distance=9mm,sibling angle=120]
			\node (A0) [draw,circle,inner sep=1] at (0:-1.5) {$\scriptstyle{1}$};
      \tikzstyle{level 1}=[counterclockwise from=0,level distance=9mm,sibling angle=120]
			\node (A1) [draw,circle,inner sep=1] at (0:0) {$\scriptstyle{2}$};
      \tikzstyle{level 1}=[counterclockwise from=0,level distance=9mm,sibling angle=120]
			\node (A2) [draw,circle,inner sep=1] at (0:1.5) {$\scriptstyle{1}$};
			\path (A0) edge [bend left=0.000000] node[el,above,font=\tiny] {$e_1$} (A1);
			\path (A1) edge [bend left=0.000000] node[el,above,font=\tiny] {$e_2$} (A2);
	\end{tikzpicture}$,
	then the stable graphs $\hat\Gamma$ obtained by splitting a vertex or adding a loop are
	$	\begin{tikzpicture}[baseline={([yshift=-.8ex]current bounding box.center)},el/.style = {inner sep=2pt, align=left, sloped},every child node/.style={inner sep=1,font=\tiny}]
      \tikzstyle{level 1}=[counterclockwise from=120,level distance=9mm,sibling angle=120]
      \node (A1) [draw,circle,inner sep=1] at (0:-1.5) {$\scriptstyle{1}$};
      \tikzstyle{level 1}=[counterclockwise from=-60,level distance=9mm,sibling angle=120]
      \node (A2) [draw,circle,inner sep=1] at (0:0) {$\scriptstyle{1}$};
      \tikzstyle{level 1}=[counterclockwise from=-60,level distance=9mm,sibling angle=120]
      \node (A3) [draw,circle,inner sep=1] at (0:1.5) {$\scriptstyle{1}$};
      \tikzstyle{level 1}=[counterclockwise from=-60,level distance=9mm,sibling angle=120]
      \node (A4) [draw,circle,inner sep=1] at (0:3) {$\scriptstyle{1}$};

			\path (A1) edge [bend left=0.000000] node[el,above,font=\tiny] {$e_1$} (A2);
			\path (A2) edge [bend left=0.000000] node[el,above,font=\tiny] {$e_3$} (A3);
			\path (A3) edge [bend left=0.000000] node[el,above,font=\tiny] {$e_2$} (A4);
	\end{tikzpicture}$,
	$	\begin{tikzpicture}[baseline={([yshift=-.4ex]current bounding box.center)},el/.style = {inner sep=2pt, align=left, sloped},every child node/.style={inner sep=1,font=\tiny}]
      \tikzstyle{level 1}=[counterclockwise from=120,level distance=9mm,sibling angle=120]
      \node (A1) [draw,circle,inner sep=1] at (0:-1.5) {$\scriptstyle{1}$};
      \tikzstyle{level 1}=[counterclockwise from=-60,level distance=9mm,sibling angle=120]
      \node (A2) [draw,circle,inner sep=1] at (0:0) {$\scriptstyle{0}$};
      \tikzstyle{level 1}=[counterclockwise from=-60,level distance=9mm,sibling angle=120]
      \node (A3) [draw,circle,inner sep=1] at (-20:1.5) {$\scriptstyle{1}$};
      \tikzstyle{level 1}=[counterclockwise from=-60,level distance=9mm,sibling angle=120]
      \node (A4) [draw,circle,inner sep=1] at (20:1.5) {$\scriptstyle{2}$};

			\path (A1) edge [bend left=0.000000] node[el,above,font=\tiny] {$e_1$} (A2);
			\path (A2) edge [bend left=0.000000] node[el,below,font=\tiny] {$e_2$} (A3);
			\path (A2) edge [bend left=0.000000] node[el,above,font=\tiny] {$e_3$} (A4);
	\end{tikzpicture}$,
	$	\begin{tikzpicture}[baseline={([yshift=-.4ex]current bounding box.center)},el/.style = {inner sep=2pt, align=left, sloped},every child node/.style={inner sep=1,font=\tiny}]
      \tikzstyle{level 1}=[counterclockwise from=120,level distance=9mm,sibling angle=120]
      \node (A1) [draw,circle,inner sep=1] at (0:-1.5) {$\scriptstyle{1}$};
      \tikzstyle{level 1}=[counterclockwise from=-60,level distance=9mm,sibling angle=120]
      \node (A2) [draw,circle,inner sep=1] at (0:0) {$\scriptstyle{1}$};
      \tikzstyle{level 1}=[counterclockwise from=-60,level distance=9mm,sibling angle=120]
      \node (A3) [draw,circle,inner sep=1] at (-20:1.5) {$\scriptstyle{1}$};
      \tikzstyle{level 1}=[counterclockwise from=-60,level distance=9mm,sibling angle=120]
      \node (A4) [draw,circle,inner sep=1] at (20:1.5) {$\scriptstyle{1}$};

			\path (A1) edge [bend left=0.000000] node[el,above,font=\tiny] {$e_1$} (A2);
			\path (A2) edge [bend left=0.000000] node[el,below,font=\tiny] {$e_2$} (A3);
			\path (A2) edge [bend left=0.000000] node[el,above,font=\tiny] {$e_3$} (A4);
	\end{tikzpicture}$,
	$\begin{tikzpicture}[baseline={([yshift=-1.9ex]current bounding box.center)},el/.style = {inner sep=2pt, align=left, sloped},every child node/.style={inner sep=1,font=\tiny}]
      \tikzstyle{level 1}=[counterclockwise from=0,level distance=9mm,sibling angle=120]
			\node (A0) [draw,circle,inner sep=1] at (0:-1.5) {$\scriptstyle{1}$};
      \tikzstyle{level 1}=[counterclockwise from=0,level distance=9mm,sibling angle=120]
			\node (A1) [draw,circle,inner sep=1] at (0:0) {$\scriptstyle{1}$};
      \tikzstyle{level 1}=[counterclockwise from=0,level distance=9mm,sibling angle=120]
			\node (A2) [draw,circle,inner sep=1] at (0:1.5) {$\scriptstyle{1}$};
			\path (A0) edge [bend left=0.000000] node[el,above,font=\tiny] {$e_1$} (A1);
			\path (A1) edge [bend left=0.000000] node[el,above,font=\tiny] {$e_2$} (A2);
			\path[every loop/.style={-}] (A1) edge [el,loop above,font=\tiny] node {$e_3$}(A1);
	\end{tikzpicture}$,
	$\begin{tikzpicture}[baseline={([yshift=-1.9ex]current bounding box.center)},el/.style = {inner sep=2pt, align=left, sloped},every child node/.style={inner sep=1,font=\tiny}]
      \tikzstyle{level 1}=[counterclockwise from=0,level distance=9mm,sibling angle=120]
			\node (A0) [draw,circle,inner sep=1] at (0:-1.5) {$\scriptstyle{0}$};
      \tikzstyle{level 1}=[counterclockwise from=0,level distance=9mm,sibling angle=120]
			\node (A1) [draw,circle,inner sep=1] at (0:0) {$\scriptstyle{2}$};
      \tikzstyle{level 1}=[counterclockwise from=0,level distance=9mm,sibling angle=120]
			\node (A2) [draw,circle,inner sep=1] at (0:1.5) {$\scriptstyle{1}$};
			\path (A0) edge [bend left=0.000000] node[el,above,font=\tiny] {$e_1$} (A1);
			\path (A1) edge [bend left=0.000000] node[el,above,font=\tiny] {$e_2$} (A2);
			\path[every loop/.style={-}] (A0) edge [el,loop above,font=\tiny] node {$e_3$}(A0);
	\end{tikzpicture}$, and
	$\begin{tikzpicture}[baseline={([yshift=-1.9ex]current bounding box.center)},el/.style = {inner sep=2pt, align=left, sloped},every child node/.style={inner sep=1,font=\tiny}]
      \tikzstyle{level 1}=[counterclockwise from=0,level distance=9mm,sibling angle=120]
			\node (A0) [draw,circle,inner sep=1] at (0:-1.5) {$\scriptstyle{0}$};
      \tikzstyle{level 1}=[counterclockwise from=0,level distance=9mm,sibling angle=120]
			\node (A1) [draw,circle,inner sep=1] at (0:0) {$\scriptstyle{2}$};
      \tikzstyle{level 1}=[counterclockwise from=0,level distance=9mm,sibling angle=120]
			\node (A2) [draw,circle,inner sep=1] at (0:1.5) {$\scriptstyle{1}$};
			\path (A0) edge [bend left=0.000000] node[el,above,font=\tiny] {$e_2$} (A1);
			\path (A1) edge [bend left=0.000000] node[el,above,font=\tiny] {$e_1$} (A2);
			\path[every loop/.style={-}] (A0) edge [el,loop above,font=\tiny] node {$e_3$}(A0);
	\end{tikzpicture}$.
\end{example}

\subsection{Compact type and rational tails} \label{sec:ct-rt2}
Just like in the pushforward case, we can calculate Deligne's pullback spectral sequence for moduli spaces of curves of compact type (respectively with rational tails) by restricting the graphs we sum over.
Unlike in the pushforward case, we now also need to restrict the differential map.
This is done by not allowing the splitting of a vertex in a way that would result in a graph with a bridge (respectively rational tail).

\section{Tautological cohomology} \label{sec:taut}
Our computations will be done using tautological cohomology, which we briefly describe here.
For any stable graph $\Gamma$ and monomial $\theta$ in $\kappa$- and $\psi$-classes pulled back from the projection onto the factors of $\Mbar_\Gamma$, we get a class in the Chow ring $A^*(\Mbar_{\Gamma})$ called a decorated stratum class. 
Let $S_{g,n}^*$ denote the $\mathbb{Q}$-vector space with the decorated stratum classes as a basis. To any decorated stratum, corresponding to $\Gamma,\theta$, we get an element,  
\[\frac{1}{|\text{Aut}(\Gamma)|}{\xi_\Gamma}_*(\theta)\in A^{*}(\Mbar_{g,n}).\]
The image of $S_{g,n}^*$ in $A^{*}(\Mbar_{g,n})$ is the tautological ring $R^*(\Mbar_{g,n})$. 
See for instance~\cite{graberpandharipande} for more details.

The image of $R^*(\Mbar_{g,n})$ in cohomology $H^{2*}(\Mbar_{g,n})$, under the cycle class map, will be denoted by $RH^*(\Mbar_{g,n})$ and be called the tautological cohomology. If $RH^r(\Mbar_{g,n})=H^{2r}(\Mbar_{g,n})$ then we will say that the cohomology (group) is tautological. 
In \cite{pixtonconjrels,PPZ}, a set $P_{g,n}^*$ of relations among the generators $S_{g,n}^*$ mapped inside $RH^*(\Mbar_{g,n})$, called Pixton's relations (and also called the generalized Faber-Zagier relations), were given and were conjectured to be all relations. 

\section{The computer implementation} \label{sec:implementation}
We have implemented Deligne's pushforward and pullback weight spectral sequences described in sections~\ref{sec:deligne} and \ref{sec:gk} into a program that takes any pair $(g,n)$ such that $2g+n \leq 12$ and computes (under the assumption that Pixton's relations are all relations, see below) the weight graded pieces of the rational cohomology of $\M_{g,n}$. The implementation is based upon knowledge of the tautological cohomology of moduli spaces of stable pointed curves and builds upon the SageMath module {\tt admcycles}, see \cite{admcycles}.

Since our implementation is based upon computations in the tautological cohomology ring, we need to make sure that the $E^1$-page of the spectral sequences consists of tautological cohomology. But if $2\tilde g+\tilde n \leq 12$ then by \cite[Theorem 1.4]{canninglarson} the cohomology of $\Mbar_{\tilde g,\tilde n}$ is tautological, and hence for all $g,n,p,q$ such that $2g+n \leq 12$, $E_1^{-p+i,q}$ and $\tilde E_1^{p+i,q}$ will be tautological. Note that since tautological cohomology is of Tate type, it follows that the cohomology of $\M_{g,n}$ will be of Tate type. Note also that the assumption that all of the cohomology of $\Mbar_{g,n}$ is tautological implies that the tautological cohomology is Gorenstein, because of Poincar\'e duality. 

In all cases where we use the Pixton relations, we have verified that these are all relations, more on this in Section~\ref{sec:pixtonrels}.

\subsection{The algorithm}
The general outline of our algorithm is as follows. Note that this does not include the refinements discussed in Section \ref{sec:stable_curves}.

\begin{itemize}
\item
As input to the algorithm, we decide which weights we compute using the pushforward spectral sequence and for which weights we instead use the pullback spectral sequence.
\item
Then, for all $\lambda\vdash n$ and all relevant graphs $\Gamma$, we compute a basis for each tautological cohomology group $(S_{g',n'}^{r'}/P^{r'}_{g',n'})^{\mathbb{S}_{\lambda'}}$ that appears in the terms of $\Bigl(\bigl(H^{q}(\Mbar_\Gamma) \otimes \mathrm{det}(E_\Gamma)\bigr)^{\mathrm{Aut}(\Gamma)}\Bigr)^{\mathbb{S}_\lambda}$.
We also keep track of how to map any element of $(S_{g',n'}^{r'})^{\mathbb{S}_{\lambda'}}$ to this basis.
Here $\mathbb{S}_\lambda:=\prod\mathbb{S}_{\lambda_i}$ is the subgroup of $\mathbb{S}_n$, permuting the marked points.
\item
The group $\mathbb{S}_\lambda$ also acts on the graphs appearing in the terms of the spectral sequences.
We choose a representative $\Gamma$ for each orbit. 
\item
Then for each such graph $\Gamma$ we compute a basis for each $\Bigl(\bigl(H^{q}(\Mbar_\Gamma) \otimes \mathrm{det}(E_\Gamma)\bigr)^{\mathrm{Aut}(\Gamma)}\Bigr)^{\mathbb{S}_\lambda}$ that appears in the spectral sequence.
This includes keeping track of how elements of $\bigl(H^{q}(\Mbar_\Gamma)\bigr)^{\mathbb S_\lambda}$ map to this basis.
\item
In the case where the $\mathbb S_\lambda$-action is trivial (i.e. $\lambda = [1,\ldots,1]$), we now compute the image of the differential.
For all other $\lambda$, we symmetrize this image. That is, we send every element of the image to a fixed representative of its $\mathbb S_\lambda$-orbit.
\item
From the ranks of the images and of the $\Bigl(\bigl(H^{q}(\Mbar_\Gamma) \otimes \mathrm{det}(E_\Gamma)\bigr)^{\mathrm{Aut}(\Gamma)}\Bigr)^{\mathbb{S}_\lambda}$, we compute the ranks of the groups on the second page of the spectral sequence up to $\mathbb{S}_\lambda$-action.
\item
We express the $\mathbb{S}_\lambda$ in terms of the irreducible representations of $\mathbb S_n$.
\end{itemize}

\begin{remark}
\label{rem:skip}
The order in which we do the steps is not exactly as described above.
In particular, we first do all computations for $\lambda = [1,\ldots,1]$, where the action is trivial.
Then, for other $\lambda$, we can restrict our computations to those that contribute to parts of the second page that are nonzero for the trivial action.
\end{remark}

During the computations we save partial results, such as the image of the page one differential restricted to a specific graph.
Thus, computing the weight graded pieces of the cohomology of $\M^\mathrm{ct}_{g,n}$ and $\M^\mathrm{rt}_{g,n}$ simply becomes a matter of combining these partial results in a different manner.
Note that this might require us to now compute partial results that we previously skipped due to Remark~\ref{rem:skip}.

The SageMath module {\tt admcycles} could already compute the generators, relations, and bases for $(S_{g,n}^r/P^r_{g,n})^{\mathbb{S}_\lambda}$.
But we added some optimizations to the original base computation code in order to reach the cases we need.
Furthermore, {\tt admcycles} can compute pushforwards and pullbacks of elements of $S_{g,n}^r/P^r_{g,n}$ under gluing maps.

\subsection{Computing the tautological ring of the moduli space of stable curves} \label{sec:stable_curves}

The main bottleneck in our computations is that of $RH^r(\Mbar_{g,n})$ appearing in the factors of $RH^*(\Mbar_\Gamma)$.
We need a basis of $RH^r(\Mbar_{g,n})$ and a way to convert any decorated stratum class to this basis.
The numbers of generators and relations grow quickly with $g, n, r$, and so solving the system of equations demands large amounts of memory already for small values of $g, n, r$.

In particular, the number of relations obtained from Pixton's formula grows similarly to the number of partitions of size at most $3r - g - 1$ (see \cite{pixtonconjrels}).
Using \cite[Proposition 2]{pixtonconjrels} one can alternatively compute the relations using the recursive structure of the boundary.
This was originally implemented by Pixton for the $\mathbb S_n$-invariant part of the cohomology.
In \cite{wenninkreconstruction} it was optimized and extended to the $\mathbb S_{n'}$-invariant part of the cohomology for $0\leq n' \leq n$.
It provides a significant improvement over directly using Pixton's formula except for very small $r$.
However it has not been implemented for general $\mathbb S_\lambda$ invariants.

Another way to deal with the large number of relations is that, if $2r > 3g-3+n$, then we can use the pairing in cases that it is perfect (which holds if $2g+n \leq 12$), to compute $RH^r(\Mbar_{g,n})$ from $RH^{3g-3+n-r}(\Mbar_{g,n})$.
This requires computing the multiplication matrix for multiplying the elements of $S_{g,n}^r$ with a basis of $RH^{3g-3+n-r}(\Mbar_{g,n})$.
The number of elements in this basis is a lot smaller than the number of Pixton's relations.

In some cases we can use the symmetry of $\Gamma$ to lower the number of decorated stratum classes we need to consider.
Let $v\in V_\Gamma$ and let $G_v\subset\mathrm{Aut}(\Gamma)$ be the subgroup generated by the automorphisms that swap the two half-edges of a self-edge of $v$ and the automorphisms that swap two self-edges of $v$.
We define
\[S^r(\Gamma, v)=(S^r_{g(v),n(v)}\otimes\mathrm{det}(E_\Gamma(v,v)))^{G_v},\]
and
\[H^*(\Gamma, v) = (H^*(\Mbar_{g(v),n(v)})\otimes\mathrm{det}(E_\Gamma(v,v)))^{G_v},\]
where $E_\Gamma(v,v)\subset E_\Gamma$ is the set of self-edges at $v$.
Because 
\[H^*(\Mbar_\Gamma,\mathrm{det}(E_\Gamma))^{\mathrm{Aut}(\Gamma)} \cong
	\Bigl(\bigl(\prod_{v\in{V_\Gamma}}H^*(\Gamma, v)\bigr)\otimes\mathrm{det}(E_\Gamma)\Bigr)^{\mathrm{Aut}(\Gamma)/(\bigoplus_v G_v)},\]
we can compute Pixton's relations on $S^r(\Gamma, v)$, rather than $S^r_{g(v),n(v)}$.

Let $m$ be the number of self-edges at $v$ and consider the subgroup $\mathbb S_2^m\subset G_v$ that permutes the half-edges of each self-edge.
Let $x\in RH^r(\Mbar_{g,n})$ and $y\in RH^{3g-3+n-r}(\Mbar_{g,n})$,
then
\[x\sum_{g\in \mathbb S_2^m}g(y) = \sum_{g\in \mathbb S_2^m}g(xy).\]
So when we use the perfect pairing, we can calculate the multiplication matrix for multiplying with $(S^r_{g(v),n(v)})^{\mathbb S_2^m}$ rather than $S^r_{g(v),n(v)}$.
Unfortunately the same does not hold for the action by the whole of $G_v$, so we cannot use it in combination with the perfect pairing.

Let $\Gamma_g$ be the graph consisting of a single vertex of genus $g$ and $5-g$ self-edges.
Here is some data that illustrates how the number of decorated stratum classes grows.
\[\begin{tabular}{lcccccc} g & $0$ & $1$ & $2$ & $3$ & $4$ & $5$ \\\midrule
	$|S^1_{g,10-2g}|$ & $512$ & $257$ & $97$ & $33$ & $11$ & $4$  \\
	$|S^2_{g,10-2g}|$ & $29836$ & $10353$ & $2404$ & $478$ & $95$ & $23$  \\
	$|(S^2_{g,10-2g})^{\mathbb S_2^{5-g}}|$ & $3959$ & $2210$ & $750$ & $229$ & $68$ & $23$  \\
	$|(S^3_{g,10-2g})^{\mathbb S_2^{5-g}}|$ & $65380$ & $29749$ & $8398$ & $1974$ & $456$ & $109$  \\
	$|(S^4_{g,10-2g})^{\mathbb S_2^{5-g}}|$ & $651710$ & $270478$ & $66650$ & $13475$ & $2581$ & $518$  \\
	$|S^4(\Gamma_g, 1)|$ & $2677$ & $6759$ & $8151$ & $5943$ & $2581$ & $518$  \\
\end{tabular}\] 
Note that the actions by $\mathbb S_2^{5-g}$ and $G_v$ provide a significant improvement.
Because of this, computing the weight graded pieces of the cohomology of $\mathcal{M}_{4,2}$ is in some respects harder than computing those of $\mathcal{M}_5$, despite the fact that the former space is of lower dimension.

\subsection{Completeness of Pixton's relations}\label{sec:pixtonrels}
In genus 0 and~1, the tautological cohomology of $\Mbar_{g,n}$ is known and Pixton's relations are all relations.
For the other $(S_{g,n}^r/P_{g,n}^r)^{\mathbb{S}_{\lambda}}$ we use, we check that
$\dim H^{2r}(\Mbar_{g,n})^{\mathbb{S}_{\lambda}} = \dim \;(S_{g,n}^r/P_{g,n}^r)^{\mathbb{S}_{\lambda}}$.
Since all cohomology is tautological this implies that Pixton's relations are all relations.
All the Betti numbers necessary for the check  
were already known except the ones for $\Mbar_{5}$.
The computation of the Betti numbers of $\Mbar_{5}$ will be part of forthcoming work by the first author together with Samir Canning, Dan Petersen and Johannes Schmitt.

In the cases where we use $H^*(\Gamma, v)$ we work with invariant parts of cohomology under a wraith group action
$(S_{g,n}^r/P_{g,n}^r)^{\mathbb S_2 \wr\mathbb S_2^{k}}$.
In the cases where $(S_{g,n}^r/P_{g,n}^r)^{\mathbb S_2^{k}}$ is computable we can check the Pixton conjecture there.
This leaves the following cases:
For $(S_{2,6}^r/P_{2,6}^r)^{\mathbb S_2^{3}}$, with $r\leq 4$ and $(S_{3,4}^5/P_{3,4}^5)^{\mathbb S_2^{2}}$, the intersection matrix over a finite field was computed.
The lower bound obtained from Pixton's relations in this way was then compared against the corresponding known Betti number.

We sum up our checks in the following proposition.
\begin{prop} \label{prop:pixton}
    For all $g$ and $n\leq N$, where $(g, N)
    \in\{(2,4),(3,3),(4,1),(5,0)\}$,
    Pixton's relations are all relations in $RH^*(\Mbar_{g,n})$.  
    Furthermore Pixton's relations are all relations in $RH^*(\Mbar_{2,5})^{\mathbb S_2}$, $RH^*(\Mbar_{2,6})^{\mathbb S_2^3}$, $RH^*(\Mbar_{3,4})^{\mathbb S_2^2}$, and $RH^*(\Mbar_{4,2})^{\mathbb S_2}$.
\end{prop}

\subsection{Bookkeeping challenges}
Implementing a simple version of the program was relatively straightforward.
The vast majority of the work went into refining it to a level where it can actually compute something that was not already known.

For the group actions we consider, whenever possible we want to work with a single representative of each orbit, rather than doing computations on the whole weighted sum of the elements of the orbit. This then results in a convention on how ordering relates to the group action:
Consider for example the action by $\mathbb{S}_2$ on $S_{g,n}^r$. If we have a representative of an orbit, that gives us enough information to reconstruct the weighted sum provided we know which two marked points are being permuted by the $\mathbb{S}_2$-action. 
We are only going to precompute a basis for $(S_{g,n}^r/P^r_{g,n})^{\mathbb{S}_2}$ once, namely with the convention that the first two points are the ones being acted upon. Thus, whenever we permute the marked points, for example when applying the permutation induced by $\phi$ from Section~\ref{sec:deligne}, we need to make sure that this convention is being respected.

Besides the action of $\mathbb{S}_\lambda$ on $H^*(\Mbar_\Gamma)$, there are the actions by $\mathrm{Aut}(\Gamma)$ and its subgroups $\mathbb S_{2}^k$ and $G_v$.
Furthermore there is the action of $\mathbb S_\lambda$ on the collection of graphs.
The bookkeeping required to make sure all the corresponding ordering conventions are respected and interact well with each other turned out to be one of the main challenges of the implementation.

While {\tt admcycles} could already compute $(S_{g,n}^r/P^r_{g,n})^{\mathbb{S}_\lambda}$, the majority of {\tt admcycles} is agnostic with respect to the $\mathbb{S}_\lambda$ action.
Similarly while the data structure of a graph in {\tt admcycles} naturally comes with an ordering of edges, as a whole {\tt admcycles} is agnostic with respect to this ordering.

\subsection{Computations}\label{sec:computations}
Using our implementation we computed the weight graded pieces of the cohomology of $\M_{g,n}$ for $n\leq N$ in the cases $(g, N) \in\{(1,6),(2,4),(3,3),(4,1),(5,0)\}$.  
We used a computer with a $32$-core $64$-thread CPU and half a terabyte of memory. The results are found in Section~\ref{sec:hodgeserre}. 

In order to compute the weight graded pieces of the cohomology of $\M_5$, we used Deligne's pushforward weight spectral sequence for $q\leq 7$ and Deligne's pullback weight spectral sequence for $q\leq 4$.
We first computed bases for all the required $RH^r(\Mbar_{g,n})$.
Here in many cases we only computed a basis for invariants under some group action as described in Section~\ref{sec:stable_curves}.
Each case could be computed within two weeks and the largest cases used almost all of the memory.
The rest of the computation was split into many smaller tasks.
The dependency graph for these tasks is not a tree. 
We wrote a simple parallel programming task scheduler to deal with this.
It also reschedules failed tasks (due to memory overload) and dynamically adjusts the number of allowed tasks when restarting the ones that failed. 
This is needed because the tasks vary a lot in terms of how much memory they require.
Furthermore, running tasks in their own processes guarantees that the memory gets released afterwards, which is relevant because {\tt admcycles} has memory leaks.
Running all these tasks took roughly 2 weeks.
Computing the weight graded pieces of the cohomology of $\M_{3,3}$ took a similar amount of time.

After computing those of $\M_{g,n}$, computing the weight graded pieces of the cohomology of $\M^\mathrm{ct}_{g,n}$ and $\M^\mathrm{rt}_{g,n}$ required relatively limited extra computation time.

In principle our program can compute all the weight graded pieces of the cohomology of $\M_{g,n}$ for all $2g+n \leq 12$ (as usual, assuming that Pixton's relations are all relations), but in practice we cannot due to memory constraints.

\section{Future directions}
In this section we will describe some cases for which one could apply (a generalization of) our program to compute weight graded pieces of the cohomology of $\M_{g,n}$.

By \cite[Remark~3.7]{canninglarsonpaynewillwacherPol}, for all $g,n$ such that $3g+2n<25$ we can restrict to the tautological contributions to the $\tilde E_1$-page in Deligne's pullback weight spectral sequence to get the correct $\tilde E_2$-page. This extends the range $2g+n \leq 12$, and it is conjectured in \cite{canninglarsonpaynewillwacherPol} to be the optimal range (this conjecture is also proven for $g+n \leq 150$). 

Assuming that Pixton's relations are all relations, we could also restrict to any single choices of $g,n,p,q$, such that all cohomology groups appearing in $E_1^{-p+i,q}$ (or $\tilde E_1^{p+i,q}$), for $i \in \{-1,0,1\}$, are tautological. This holds for instance for $\tilde E_1^{p+i,q}$ for any $g,n,p$ and even $q \leq 4$, and also for $q=6$ when assuming that $g \geq 10$, see \cite[Theorem 1.5]{canninglarsonpayneExt}. See further in \cite[Theorem~8]{canningGor} for some other cases with $q$ even and $g \leq 8$. 

Presently, the program only works with tautological cohomology, but it would in principle be generalizable if one knows how the non-tautological classes behave under gluing maps. An example would be $H^{11}(\Mbar_{g,n})$, which is non-zero only for $g=1$ and $n \geq 11$ (see \cite[Theorem~1.1]{CLP}), and for the non-zero classes the pullback under gluing maps is described in \cite[Section 3]{paynewillwacherweight11}. Using this description a possible extension of our program would compute the weight $11$ graded piece of the compactly supported cohomology of $\M_{g,n}$ for any $g,n$ (since all odd contributions up to weight $9$ to the cohomology of $\Mbar_{\tilde g,\tilde n}$ for any $\tilde g,\tilde n$ vanish, see \cite{BFP} and \cite{paynewillwacherweight11}).

\section{A list of $\mathbb S_n$-equivariant Hodge-Serre polynomials} \label{sec:hodgeserre}
Here we collect the results from our computations described in Section~\ref{sec:computations}. 
This notation will be slightly different from the one in Theorems~\ref{thm:main1} and \ref{thm:main2}:
For any $g,n \geq 0$, with $2g-2+n > 0$, define the polynomial  
\[
p_{g,n}=\sum_{r \geq 0} \sum_{q \geq 0} \sum_{\lambda \vdash n} \dim \left(U^{2q,r,\lambda}_{g,n} \right)t^r \L^q s_{\lambda}, 
\]
with the notation from Section~\ref{sec:introduction}. We will only compute this polynomial in cases for which all of the weight graded pieces of the cohomology of $\M_{g,n}$ is Tate, which justifies the use of $\L$ as variable and that only even weights are considered. 

\begin{thm} We have that,
\begin{align*}
p_{1,1} &= s_{1},\\ 
p_{1,2} &= s_{2},\\ 
p_{1,4} &= t^{3} \L^{3} s_{2,1,1} + t^{4} \L^{4} s_{3,1} + \left(t^{3} \L^{2} + 1\right) s_{4} 
,\\ 
p_{1,5} &= t^{5} \L^{5} s_{1,1,1,1,1} + t^{5} \L^{5} s_{2,2,1} + \left(t^{4} \L^{4} + t^{3} \L^{3}\right) s_{3,1,1} + \left(t^{5} \L^{5} + t^{4} \L^{4} + t^{4} \L^{3}\right) s_{3,2} + \left(t^{4} \L^{4} + t^{3} \L^{2}\right) s_{4,1} \\
&\quad+ \left(t^{5} \L^{5} + t^{3} \L^{2} + 1\right) s_{5}
,\\ 
p_{1,6} &= t^{5} \L^{5} s_{2,1,1,1,1} + t^{5} \L^{5} s_{2,2,1,1} + \left(t^{6} \L^{6} + t^{5} \L^{5}\right) s_{2,2,2} + 2 t^{6} \L^{6} s_{3,1,1,1} + \left(t^{6} \L^{6} + 2 t^{5} \L^{5} + t^{5} \L^{4} + t^{4} \L^{4} + t^{4} \L^{3}\right) s_{3,2,1} \\
&\quad+ \left(t^{5} \L^{5} + t^{4} \L^{4} + t^{4} \L^{3}\right) s_{3,3} + \left(t^{6} \L^{6} + t^{5} \L^{4} + t^{4} \L^{4} + t^{3} \L^{3}\right) s_{4,1,1} + \left(t^{6} \L^{6} + 2 t^{5} \L^{5} + t^{4} \L^{4} + t^{4} \L^{3} + t^{3} \L^{2}\right) s_{4,2} \\
&\quad+ \left(t^{5} \L^{5} + t^{4} \L^{4} + t^{4} \L^{3} + t^{3} \L^{2}\right) s_{5,1} + \left(t^{5} \L^{5} + t^{3} \L^{2} + 1\right) s_{6},\\ 
p_{2,0} &= 1 
,\\ p_{2,1} &= \left(t^{2} \L + 1\right) s_{1} 
,\\ p_{2,2} &= t^{2} \L s_{1,1} + \left(t^{5} \L^{5} + t^{2} \L + 1\right) s_{2} 
,\\ p_{2,3} &= \left(t^{5} \L^{5} + t^{2} \L\right) s_{2,1} + \left(t^{5} \L^{5} + t^{3} \L^{2} + t^{2} \L + 1\right) s_{3} 
,\\ p_{2,4} &= \left(t^{7} \L^{7} + t^{6} \L^{5}\right) s_{2,1,1} + t^{5} \L^{5} s_{2,2} + \left(t^{5} \L^{5} + t^{3} \L^{2} + t^{2} \L\right) s_{3,1} + \left(t^{8} \L^{7} + t^{6} \L^{6} + t^{5} \L^{5} + t^{5} \L^{3} + t^{3} \L^{2} + t^{2} \L + 1\right) s_{4}
,\\ p_{3,0} &= t^{6} \L^{6} + t^{2} \L + 1 
,\\ p_{3,1} &= \left(t^{8} \L^{7} + t^{6} \L^{6} + t^{4} \L^{2} + 2 t^{2} \L + 1\right) s_{1} 
,\\ p_{3,2} &= \left(t^{8} \L^{7} + t^{4} \L^{2} + t^{2} \L\right) s_{1,1} + \left(t^{8} \L^{7} + t^{6} \L^{6} + t^{5} \L^{3} + t^{4} \L^{2} + 2 t^{2} \L + 1\right) s_{2} 
,\\ p_{3,3} &= \left(2 t^{8} \L^{7} + 2 t^{5} \L^{3} + t^{4} \L^{2} + t^{2} \L\right) s_{2,1} + \left(t^{9} \L^{9} + 2 t^{9} \L^{8} + 2 t^{8} \L^{7} + t^{6} \L^{6} + 2 t^{5} \L^{3} + t^{4} \L^{2} + 2 t^{2} \L + 1\right) s_{3} 
,\\ p_{4,0} &= t^{5} \L^{3} + t^{4} \L^{2} + t^{2} \L + 1 
,\\ p_{4,1} &= \left(t^{9} \L^{8} + t^{7} \L^{4} + t^{6} \L^{3} + 2 t^{5} \L^{3} + 2 t^{4} \L^{2} + 2 t^{2} \L + 1\right) s_{1} 
,\\ p_{5,0} &= t^{14} \L^{12} + t^{7} \L^{4} + t^{6} \L^{3} + t^{5} \L^{3} + t^{4} \L^{2} + t^{2} \L + 1 
.
\end{align*}
\end{thm}

\printbibliography
\end{document}